\documentclass[11pt,leqno]{article}
\usepackage{graphicx, amsfonts, amsthm, amsxtra, amssymb, verbatim}
\usepackage[mathscr]{euscript}

\begin{document}
\def\Sp{{\rm Sp}}  
\def\Ra{{\rm Ra}}
\def\GM{{\rm GM}}

\def\per{{\sf pm}}      
\def\perr{{\sf q}}        
\def\perdo{{\cal K}}   
\def\sfl{{\mathrm F}} 
\def\sp{{\mathbb S}}  
 
\newcommand\diff[1]{\frac{d #1}{dz}} 
\def\End{{\rm End}}              

\def\sing{{\rm Sing}}            
\def\spec{{\rm Spec}}            
\def\cha{{\rm char}}             
\def\Gal{{\rm Gal}}              
\def\jacob{{\rm jacob}}          
\def\tjurina{{\rm tjurina}}      
\newcommand\Pn[1]{\mathbb{P}^{#1}}   
\def\Ff{\mathbb{F}}                  
\def\Z{\mathbb{Z}}                   
\def\Gm{\mathbb{G}_m}                 
\def\Q{\mathbb{Q}}                   
\def\C{\mathbb{C}}                   
\def\O{{\cal O}}                     
\def\as{\mathbb{U}}                  
\def\ring{{\mathsf R}}                         
\def\R{\mathbb{R}}                   
\def\N{\mathbb{N}}                   
\def\A{\mathbb{A}}                   
\def\uhp{{\mathbb H}}                
\newcommand\ep[1]{e^{\frac{2\pi i}{#1}}}
\newcommand\HH[2]{H^{#2}(#1)}        
\def\Mat{{\rm Mat}}              
\newcommand{\mat}[4]{
     \begin{pmatrix}
            #1 & #2 \\
            #3 & #4
       \end{pmatrix}
    }                                
\newcommand{\matt}[2]{
     \begin{pmatrix}                 
            #1   \\
            #2
       \end{pmatrix}
    }
\def\ker{{\rm ker}}              
\def\cl{{\rm cl}}                
\def\dR{{\rm dR}}                

\def\hc{{\mathsf H}}                 
\def\Hb{{\cal H}}                    
\def\GL{{\rm GL}}                
\def\pese{{\sf P}}                  
\def\pedo{{\cal  P}}                  
\def\PP{\tilde{\cal P}}              
\def\cm {{\cal C}}                   
\def\K{{\mathbb K}}                  
\def\k{{\mathsf k}}                  
\def\F{{\cal F}}                     
\def\M{{\cal M}}
\def\RR{{\cal R}}
\newcommand\Hi[1]{\mathbb{P}^{#1}_\infty}
\def\pt{\mathbb{C}[t]}               
\def\W{{\cal W}}                     
\def\gr{{\rm Gr}}                
\def\Im{{\rm Im}}                
\def\Re{{\rm Re}}                
\def\depth{{\rm depth}}
\newcommand\SL[2]{{\rm SL}(#1, #2)}    
\newcommand\PSL[2]{{\rm PSL}(#1, #2)}  
\def\Resi{{\rm Resi}}              

\def\L{{\cal L}}                     
\def\Aut{{\rm Aut}}              
\def\any{R}                          
\newcommand\ovl[1]{\overline{#1}}    

\def\T{{\cal T }}                    
\def\tr{{\mathsf t}}                 
\newcommand\mf[2]{{M}^{#1}_{#2}}     
\newcommand\bn[2]{\binom{#1}{#2}}    
\def\ja{{\rm j}}                 
\def\Sc{\mathsf{S}}                  
\newcommand\es[1]{g_{#1}}            
\newcommand\V{{\mathsf V}}           
\newcommand\WW{{\mathsf W}}          
\newcommand\Ss{{\cal O}}             
\def\rank{{\rm rank}}                
\def\Dif{{\cal D}}                   
\def\gcd{{\rm gcd}}                  
\def\zedi{{\rm ZD}}                  
\def\BM{{\mathsf H}}                 
\def\plf{{\sf pl}}                             
\def\sgn{{\rm sgn}}                      
\def\diag{{\rm diag}}                   
\def\hodge{{\rm Hodge}}
\def\HF{{\sf F}}                                
\def\WF{{\sf W}}                               
\def\HV{{\sf HV}}                                
\def\pol{{\rm pole}}                               
\def\bafi{{\sf r}}
\def\codim{{\rm codim}}                               
\def\id{{\rm id}}                               
\def\gms{{\sf M}}                           
\def\Iso{{\rm Iso}}                           
\newtheorem{theo}{Theorem}
\newtheorem{exam}{Example}
\newtheorem{coro}{Corollary}
\newtheorem{defi}{Definition}
\newtheorem{prob}{Problem}
\newtheorem{lemm}{Lemma}
\newtheorem{prop}{Proposition}
\newtheorem{rem}{Remark}
\newtheorem{conj}{Conjecture}
\newtheorem{calc}{}

\begin{center}
{\LARGE\bf  Modular-type functions attached to mirror quintic Calabi-Yau varieties 
\footnote{ 
Math. classification: 14N35, 
14J15, 32G20
}
}
\\
\vspace{.25in} {\large {\sc Hossein Movasati}} \\
Instituto de Matem\'atica Pura e Aplicada, IMPA, \\
Estrada Dona Castorina, 110,\\
22460-320, Rio de Janeiro, RJ, Brazil, \\
{\tt www.impa.br/$\sim$ hossein, hossein@impa.br} 
\end{center}
\begin{abstract}
In this article we study a differential algebra  of modular-type functions attached to the periods of a one parameter family 
of Calabi-Yau varieties which is mirror dual to the universal family of quintic threefolds. Such an algebra is generated by seven functions
satisfying functional and differential equations in parallel to  the modular functional equations of classical Eisenstein series and the 
Ramanujan differential equation. Our result is the first
example of automorphic-type functions attached to varieties whose period domain is not Hermitian symmetric. 
It is a reformulation and realization of a problem of Griffiths around seventies on the existence of automorphic functions for the moduli of polarized Hodge structures. 
\end{abstract}
\section{Introduction}
In 1991 Candelas, de la Ossa, Green and Parkes in \cite{can91}  calculated in the framework of mirror symmetry 
a generating function, called the Yukawa coupling,  which predicts the number of rational curves of a fixed degree in a generic quintic 
threefold. Since then there was many efforts to relate the Yukawa 
coupling to the  classical modular or quasi-modular forms, however, there was no success. 
The theory of modular  or quasi-modular forms is attached to elliptic curves, or better to say their periods, see for instance \cite{kon, ho14}. In general
the available holomorphic automorphic function theories are attached to varieties whose Hodge structures form a Hermitian symmetric domain. 
This is not the case for mirror quintic Calabi-Yau threefolds which are the  underlying varieties of the Yukawa coupling.
An attempt to formulate automorphic function theories beyond Hermitian symmetric domains was first done around 
seventies by P. Griffiths  in the framework of Hodge structures, see \cite{gr70}. However, such a formulation has  lost the generating 
function role of modular forms. 
The main aim of the present text is to reformulate and realize the construction of a modular-type function theory attached to mirror quintic threefolds.
The present work continues and simplifies the previous article \cite{ho21}.

Consider the following fourth order linear differential equation:
\begin{equation}
\label{4jul2011}
\theta^4-z(\theta+\frac{1}{5})(\theta+\frac{2}{5})(\theta+\frac{3}{5})(\theta+\frac{4}{5})=0,\  \ \theta=z\frac{\partial}{\partial z}.
\end{equation}
A basis of the solution space of (\ref{4jul2011}) is given by:
$$
\psi_i(z)=\frac{1}{i!}\frac{\partial^ i}{\partial \epsilon^ i}(5^{-5\epsilon}F(\epsilon,z)),\ \ j=0,1,2,3,
$$ 
where
$$
F(\epsilon,z):=\sum_{n=0}^\infty \frac{(\frac{1}{5}+\epsilon)_n(\frac{2}{5}+\epsilon)_n(\frac{3}{5}+\epsilon)_n (\frac{4}{5}+\epsilon)_n }{(1+\epsilon)_n^4}z^{\epsilon+n}
$$
and $(a)_n:=a(a+1)\cdots(a+n-1)$ for  $n>0$ and  $(a)_0:=1$.
We use the base change
$$
\begin{pmatrix}
 x_{11}\\
x_{21}\\
x_{31}\\
x_{41}
\end{pmatrix}=
\begin{pmatrix}
   0 &0&1&0\\
0&0&0&1\\
0&5&\frac{5}{2}&-\frac{25}{12}\\
-5&0&-\frac{25}{12}&200\frac{\zeta(3)}{(2\pi i)^3}
  \end{pmatrix}\begin{pmatrix}
\frac{1}{5^4}\psi_3\\
\frac{2\pi i}{5^ 4}\psi_2\\
\frac{(2\pi i)^2}{5^ 4}\psi_1\\
\frac{(2\pi i)^ 3}{5^ 4}\psi_0
\end{pmatrix}.
$$
In the new basis $x_{1i}$ the monodromy of (\ref{4jul2011}) around the singularities 
$z=0$ and $z=1$ are respectively given by:  
$$
M_0:=
\begin{pmatrix}
1&1& 0& 0\\
0&1& 0& 0\\
5&5& 1& 0\\
0&-5&-1&1 
\end{pmatrix},\ \ \
M_1:=
\begin{pmatrix}
1&0&0&0\\
0&1&0&1\\
0&0&1&0\\
0&0&0&1 
\end{pmatrix}
$$
that is, the analytic continuation of the $4\times 1$ matrix $[x_{i1}]$ around the singularity $z=0$ (respectively $z=1$) is given by $M_0[x_{i1}]$, respectively $M_1[x_{i1}]$ 
(see for instance \cite{dormor}, \cite{vanvan} and \cite{ChYaYu} for similar calculations). The functions $x_{i1}$ can be written as 
periods of a holomorphic differential 3-forms over topological cycles with integral coefficients (see \S\ref{section3}).
 
Let 
$$
\tau_0:=\frac{x_{11}}{x_{21}},\ \  q:=e^{2\pi i \tau_0},
$$
and 
$$
x_{ij}:=\theta^{j-1}x_{i1},\ \ i,j=1,2,3,4.
$$
\begin{theo}
\label{main}
Let
\begin{eqnarray*}
t_0 &=&  x_{21},
\\
t_1 &=&5^4x_{21}\left((6z-1)x_{21}+5(11z-1)x_{22}+25(6z-1)x_{23}+125(z-1)x_{24}\right ),
\\
t_2 &=& 5^4x_{21}^2 \left ((2z-7)x_{21}+15(z-1)x_{22}+25(z-1)x_{23}\right ),
\\
t_3 &=&  5^4x_{21}^3\left (
   (z-6)x_{21}+5(z-1)x_{22} \right ),
\\
t_4 &=&
   zx_{21}^5,
\\
t_5 &=&5^5(z-1)x_{21}^2\left(x_{12}x_{21}-x_{11}x_{22} \right ),
\\
t_6 &=& 5^5(z-1)x_{21}\left ( 
  3(x_{12}x_{21}-x_{11}x_{22})+ 5(x_{13}x_{21}-x_{11}x_{23}) \right ).
\end{eqnarray*}
There are holomorphic functions $h_i$ defined in some neighborhood of $0\in\C$ such that
\begin{equation}
 \label{18J2011}
t_i=(\frac{2\pi i}{5})^{d_i}h_i(e^{2\pi i \tau_0}),
\end{equation}
where
$$
d_i:=3(i+1),\ i=0,1,2,3,4, \ \ d_5:=11,\ \  d_6:=8.
$$
Moreover, $t_i$'s satisfy the following ordinary differential equation:
\begin{equation}
\Ra:\ \
\label{lovely}
\left \{ \begin{array}{l}
\dot t_0=\frac{1}{t_5}
(6\cdot 5^4t_0^5+t_0t_3-5^4t_4)
\\
\dot t_1=\frac{1}{t_5}
(-5^8t_0^6+5^5t_0^4t_1+5^8t_0t_4+t_1t_3)
 
\\
\dot t_2=\frac{1}{t_5}
(-3\cdot 5^9t_0^7-5^4t_0^5t_1+2\cdot 5^5t_0^4t_2+3\cdot 5^9 t_0^2t_4+5^4t_1t_4+2t_2t_3)
\\
\dot t_3=\frac{1}{t_5}
(-5^{10}t_0^8-5^4t_0^5t_2+3\cdot 5^5t_0^4t_3+5^{10}t_0^3t_4+5^4t_2t_4+3t_3^2)

\\
\dot t_4=\frac{1}{t_5}
(5^6t_0^4t_4+5t_3t_4)
\\
\dot t_5=\frac{1}{t_5}
(-5^4t_0^5t_6+3\cdot 5^5t_0^4t_5+2t_3t_5+5^4t_4t_6)
\\
\dot t_6=\frac{1}{t_5}
(3\cdot 5^5t_0^4t_6-5^5t_0^3t_5-2t_2t_5+3t_3t_6)
  
\end{array} \right.
\end{equation}
with $\dot *:=\frac{\partial *}{\partial\tau_0}$.
%
%
%
\end{theo}
We define $\deg(t_i)=d_i$. In this way the  right hand side of $\Ra$ we have homogeneous rational functions of degree 
$4,7,10,13,16,12,9$ which 
is compatible with the left hand side if we assume that the derivation 
increases the degree by one. The ordinary differential equation $\Ra$ is a generalization of the Ramanujan differential 
equation between Eisenstein series, see for instance \cite{zag123, ho14}.

We write each $h_i$ as a formal power series in $q$, $h_i=\sum_{n=0}^\infty h_{i,n}q^n$ and substitute in (\ref{lovely}) with 
$\dot *=5q\frac{\partial *}{\partial q}$  and we see that 
it determines all the coefficients $h_{i,n}$ uniquely with the initial values:
\begin{equation}
 \label{22july2010}
h_{0,0}=\frac{1}{5}, \ h_{0,1}=24,\  h_{4,0}=0
\end{equation}
and assuming that $h_{5,0}\not =0$. In fact the differential equation (\ref{lovely}) seems to be the simplest way of writing  
the mixed recursion between 
$h_{i,n},\ n\geq 2$. 
Some of the first coefficients of $h_i$'s are given in the table at the end of the Introduction.
The differential Galois group of (\ref{4jul2011}) is ${\rm Sp}(4,\C)$. This together with the equality (\ref{18J2011}) imply that
the functions  $h_i,\ i=0,1,\ldots,6$ are algebraically independent over $\C$ (see \cite{ho21} Theorem 2).

The reader who is expert in classical modular forms may ask for the functional equations of $t_i$'s. 
Let $\uhp$ be the monodromy covering of $(\C-\{0,1\})\cup\{\infty\}$ associated to the monodromy group 
$\Gamma:=\langle M_1,M_0\rangle$ of (\ref{4jul2011}) (see \S \ref{18july2011}). 
The set $\uhp$ is conjecturally  biholomorphic to the upper half plane. This is equivalent to say that the only relation between 
$M_0$ and $M_1$ is $(M_0M_1)^5=I$. We do not need to assume this conjecture because we do not need 
the coordinate system on $\uhp$ given by this biholomorphism. 
The monodromy group  $\Gamma$ acts from the left on $\uhp$ in a canonical way: 
$$
(A,w)\mapsto A(w)\in\uhp,\ \ \ A\in\Gamma,\ w\in\uhp
$$
and the quotient $\Gamma\backslash \uhp$ is biholomorphic to 
$(\C-\{0,1\})\cup\{\infty\}$. This action has one elliptic point $\infty$ of order $5$ and two cusps $0$ and $1$.
We can regard $x_{ij}$ as holomorphic one valued functions on $\uhp$. For simplicity we use the same notation for these functions: 
$x_{ij}:\uhp\to \C$. We define
$$
\tau_i:\uhp\to \C,\ \ \tau_0:=\frac{x_{11}}{x_{21}},\ \tau_1:=\frac{x_{31}}{x_{21}},\  \tau_2:=\frac{x_{41}}{x_{21}},\
\tau_3:=\frac{x_{31}x_{22}-x_{32}x_{21}}{x_{11}x_{22}-x_{12}x_{21}}
$$ 
which are a priori meromorphic functions on $\uhp$. 
We will use $\tau_0$ as a local coordinate around a point $w\in\uhp$ whenever $w$ is not a pole of $\tau_0$ and
the derivative of $\tau_0$ does not vanish at $w$. In this way we need to express $\tau_i,\ i=1,2,3$ as functions of $\tau_0$:
\begin{eqnarray*}
\tau_1 &=& -\frac{25}{12}+\frac{5}{2}\tau_0(\tau_0+1)+\frac{\partial H}{\partial \tau_0}, \\ 
\tau_2 &=& 200\frac{\zeta(3)}{(2\pi i)^3}-\frac{5}{6}\tau_0(\frac{5}{2}+\tau_0^2)-\tau_0\frac{\partial H}{\partial \tau_0}-2H,\\
\tau_3 &=& \frac{\partial \tau_1}{\partial \tau_0},\\
\end{eqnarray*}
where
\begin{equation}
\label{dastgir}
H=\frac{1}{(2\pi i)^3} \sum_{n=1}^\infty (\sum_{d|n}n_d d^3)\frac{e^{2\pi i \tau_0 n}}{n^3}
\end{equation}
and $n_d$ is the virtual number of rational curves of degree $d$ in a generic quintic threefold. 
The numbers $n_d$ are also called instanton numbers or BPS degeneracies.  
A complete description of the image of $\tau_0$ is not yet known. 
Now, $t_i$'s  are well-defined holomorphic functions on $\uhp$.
 The functional equations of $t_i$'s  with respect to the action of an arbitrary element of  $\Gamma$ are complicated mixed 
equalities which we have described in \S\ref{18july2011}.
Since $\Gamma$ is generated by $M_0$ and $M_1$ it is enough to explain them for  these two elements. 
The functional equations of $t_i$'s with respect to the action of $M_0$ and written in the 
$\tau_0$-coordinate are the trivial equalities $t_i(\tau_0)=t_i(\tau_0+1),\ i=0,1,\ldots,6$. 
\begin{theo}
\label{main2}
With respect to the action of $M_1$, $t_i$'s written in the $\tau_0$-coordinate satisfy the following functional equations:  
  \begin{eqnarray*}
 t_0(\tau_0) &=& t_0(\frac{\tau_0}{\tau_2+1})\frac{1}{\tau_2+1}, \\
t_1(\tau_0) &=&   t_1(\frac{\tau_0}{\tau_2+1})\frac{1}{(\tau_2+1)^2}+t_7(\frac{\tau_0}{\tau_2+1})\frac{\tau_0\tau_3-\tau_1}{(\tau_2+1)(\tau_0^2\tau_3-\tau_0\tau_1+\tau_2+1)} +\\
& & t_{9}(\frac{\tau_0}{\tau_2+1})\frac{-\tau_0}{ (\tau_2+1)^2}+\frac{1}{\tau_2+1},\\
t_2(\tau_0)  &=&  t_2(\frac{\tau_0}{\tau_2+1})\frac{1}{(\tau_2+1)^3}+t_6(\frac{\tau_0}{\tau_2+1})\frac{\tau_0\tau_3-\tau_1}{(\tau_2+1)^2(\tau_0^2\tau_3-\tau_0\tau_1+\tau_2+1)} +\\
 & & t_8(\frac{\tau_0}{\tau_2+1})\frac{-\tau_0}{ (\tau_2+1)^3},\\
t_3(\tau_0)  &=& t_3(\frac{\tau_0}{\tau_2+1})\frac{1}{(\tau_2+1)^4}+t_5(\frac{\tau_0}{\tau_2+1})\frac{\tau_0\tau_3-\tau_1}{(\tau_2+1)^3(\tau_0^2\tau_3-\tau_0\tau_1+\tau_2+1)} \\
t_4(\tau_0)  &=& t_4(\frac{\tau_0}{\tau_2+1})\frac{1}{(\tau_2+1)^5},\\
t_5(\tau_0) &=&   t_5(\frac{\tau_0}{\tau_2+1})\frac{1}{(\tau_2+1)^2(\tau_0^2\tau_3-\tau_0\tau_1+\tau_2+1)}, \\
t_6(\tau_0) &=&  t_6(\frac{\tau_0}{\tau_2+1})\frac{1}{(\tau_2+1)(\tau_0^2\tau_3-\tau_0\tau_1+\tau_2+1)}+t_8(\frac{\tau_0}{\tau_2+1}) \frac{\tau_0^2}{(\tau_2+1)^3},\\
\end{eqnarray*}
where
\begin{eqnarray*}
t_7 &:=& \frac{(5^5t_0^4+t_3)t_6-(5^5t_0^3+t_2)t_5}{5^4(t_4-t_0^5)}, \\ 
t_8 &:=& \frac{5^4(t_0^5-t_4)}{t_5},\\
t_{9} &:=& \frac{-5^5t_0^4-t_3}{t_5}.
\end{eqnarray*}
 \end{theo}
Since we do not know about the global behavior of $\tau_0$, the above equalities must be interpreted in the following way: 
for any fixed branch of $t_i(\tau_0)$ there
is a path $\gamma$ in the image of $\tau_0: \uhp \to \C$ which connects $\tau_0$ to $\frac{\tau_0}{\tau_2+1}$ and such that 
the analytic  continuation of $t_i$'s along the path $\gamma$ satisfy the above equalities. 
For this reason it may be reasonable to use a new name for all $t_i$'s in the right hand side of the equalities in Theorem \ref{main2}.
We could also state Theorem \ref{main2} without using any local
coordinate system on $\uhp$: we regard $t_i$ as holomorphic functions on $\uhp$ and, for instance, the first equality in 
Theorem \ref{main2} can be derived from the equalities:
$$
t_0(w)=t_0(M_1(w))\frac{1}{\tau_2(w)+1},\ \ \ \tau_0(M_1(w))=\frac{\tau_0(w)}{\tau_2(w)+1}.
$$

The Yukawa coupling $Y$ turns out to be
\begin{eqnarray*}
Y &=& \frac{5^8(t_4-t_0^5)^2}{t_5^3}\\
&=& (\frac{2\pi i}{5})^{-3}\left (5+2875 \frac{q}{1-q}+ 609250\cdot 2^3\frac{q^2}{1-q^2}+
\cdots+ n_d d^3\frac{q^d}{1-q^d}+\cdots\right )
\end{eqnarray*}
and so it satisfies the functional equation
$$
Y(\tau_0)=Y( \frac{\tau_0}{\tau_2+1}) \frac{(\tau_0^2\tau_3-\tau_0\tau_1+\tau_2+1)^3}{(\tau_2+1)^4}.
$$
The basic idea behind all the computations in Theorem \ref{main} and Theorem \ref{main2} lies in the following geometric theorem:
\begin{theo}
 \label{main3}
Let $T$  be the moduli of pairs $(W,[\alpha_1,\alpha_2,\alpha_3,\alpha_4])$, where $W$ is a mirror quintic Calabi-Yau threefold and
$$
\alpha_i\in F^{4-i}\backslash F^{5-i},\ \ i=1,2,3,4,
$$
$$ 
[\langle \alpha_i,\alpha_j\rangle]=\Phi.
$$
Here, $H_\dR^3(W)$ is the third algebraic de Rham cohomology of $W$, 
$F^i$ is the $i$-th piece of the Hodge filtration of $H^3_\dR(W)$, $\langle\cdot,\cdot \rangle$ is the intersection form 
in $H^3_\dR(W)$ and $\Phi$ is 
the constant matrix:
\begin{equation}
\label{31aug10}
\Phi:=
\begin{pmatrix}
 0&0&0&1\\
0&0&1&0\\
0&-1&0&0\\
-1&0&0&0
\end{pmatrix}.
\end{equation}
Then there is a unique vector field $\Ra$ in $T$ such the Gauss-Manin connection of the universal family of mirror quintic 
Calabi-Yau varieties over $T$ composed 
with the vector field $\Ra$, namely $\nabla_\Ra$, satisfies:
\begin{eqnarray*}
\nabla_{\Ra}(\alpha_1) &=& \alpha_2,\\
\nabla_{\Ra}(\alpha_2)  &=& Y\alpha_3,\\
\nabla_{\Ra}(\alpha_3) &=& -\alpha_4,\\
 \nabla_{\Ra}(\alpha_4) &=& 0 
\end{eqnarray*}
for some regular function $Y$ in $T$. In fact,
\begin{equation}
 \label{thanksdeligne}
T\cong\{(t_0,t_1,t_2,t_3,t_4,t_5,t_6)\in \C^7\mid t_5t_4(t_4-t_0^5)\not =0\},
\end{equation} 
and under this isomorphism  the vector field $\Ra$ as an 
ordinary differential equation
is  (\ref{lovely}) and $Y=\frac{5^8(t_4-t_0)^2}{t_5^3}$ is the Yukawa coupling.
\end{theo}
The space of classical (or elliptic) modular or quasi-modular forms of a fixed degree for discrete subgroups of $\SL 2\R$, is finite. 
This simple observation is the origin of many number theoretic applications (see for instance \cite{zag123}) and we may try to generalize 
such applications to the context of present text.
Having this in mind, we have to describe the behavior of $t_i$'s around the other cusp $z=1$. This will be done in another article.
All the calculations of the present text, and in particular the calculations of $p_i$'s and the 
differential equation (\ref{lovely}), are done by Singular, see \cite{GPS01}. 
The reader who does not want to calculate everything by his 
own effort can obtain the corresponding Singular code from my web page. 
Many arguments of the the present text work for an arbitrary Calabi-Yau differential equation in the sense of \cite{alenstzu}. 
In this paper we mainly focus on the geometry of mirror quintic Calabi-Yau varieties which led us to explicit calculations. 
Therefore, the results for an arbitrary Calabi-Yau equation 
is postponed to another paper.

The present text is organized in the following way. \S\ref{section1} is dedicated to the algebro-geometric aspects of 
mirror quintic Calabi-Yau threefolds. In this section we describe how one can get the differential equation (\ref{lovely}) using
the Gauss-Manin connection of families of mirror quintic Calabi-Yau varieties enhanced with elements in their de Rham cohomologies.
Theorem \ref{main3} is proved in this section.
In \S \ref{section2} we describe a solution of (\ref{lovely}). This solution is characterized by a special format of the period 
matrix of mirror quintic Calabi-Yau varieties. Finally in \S\ref{section3} we describe such a solution in terms of the periods 
$x_{ij}$ and the corresponding  $q$-expansion. In this section we first prove Theorem \ref{main2} and then Theorem \ref{main}. 

Don Zagier  pointed out that using the parameters $t_7,t_8,t_9$ 
the differential equation (\ref{lovely}) must look simpler. I was able to rewrite it in the following way:
\begin{equation}
\label{lovely1}
\left \{ \begin{array}{l}
\dot t_0=t_8-t_0t_9\\
\dot t_1=-t_1t_9-5^4t_0t_8\\
\dot t_2=-t_1t_8-2t_2t_9-3\cdot 5^5t_0^2t_8\\
\dot t_3=4t_2t_8-3t_3t_9-5(t_7t_8-t_9t_6)t_8\\
\dot t_4=-5t_4t_9\\
\dot t_5=-t_6t_8-3t_5t_9-t_3\\
\dot t_6= -2t_6t_9-t_2-t_7t_8\\
\dot t_7=-t_7t_9-t_1\\
\dot t_8=\frac{t_8^2}{t_5}t_6-3t_8t_9\\
\dot t_9= \frac{t_8^2}{t_5} t_7-t_9^2
\end{array} \right.
\end{equation}
My sincere thanks go to Charles Doran, Stefan Reiter and Duco van Straten  for useful conversations and their
interest on the topic of the present text. In particular, I would like to thank Don Zagier whose comments on the 
first draft of the text motivated me to write Theorem \ref{main} and Theorem \ref{main2} in a 
more elementary way and without geometric considerations. I found it useful for me and a reader who seeks 
for number theoretic  applications similar to those for classical modular forms, see for instance \cite{zag123}.  
I would also like to thank both mathematics
institutes IMPA and MPIM for providing excellent research ambient during the preparation of the present text.

{
\tiny
\begin{center}
\begin{tabular}{|c|c|c|c|c|c|c|c|}
\hline
&
$q^0$
&
$q^1$
&
$q^2$
&
$q^3$
&
$q^4$
&
$q^5$
&
$q^6$
\\ \hline
$\frac{1}{24}h_0$ &
   $\frac{1}{120}$
&
   1
&
   175
&
   117625
&
   111784375
&
   126958105626
&
   160715581780591
\\ \hline
 $\frac{-1}{750}h_1$ &
   $\frac{1}{30}$
&
   3
&
   930
&
   566375
&
   526770000
&
   592132503858
&
   745012928951258
\\ \hline
$\frac{-1}{50}h_2$ &
   $\frac{7}{10}$
&
   107
&
   50390
&
   29007975
&
   26014527500
&
   28743493632402
&
   35790559257796542

\\ \hline
$\frac{-1}{5}h_3$ &
$\frac{6}{5}$
&
   71
&
   188330
&
   100324275
&
   86097977000
&
   93009679497426
&
   114266677893238146
\\ \hline
$-h_4$ &
0
&
   -1
&
   170
&
   41475
&
   32183000
&
   32678171250
&
   38612049889554
\\ \hline
$\frac{1}{125}h_5$ &
 $-\frac{1}{125}$
&
   15
&
   938
&
   587805
&
   525369650
&
   577718296190
&
   716515428667010

\\  \hline
$\frac{1}{25}h_6$ &
-$\frac{3}{5}$ &
187& 
28760 &
16677425 & 
15028305250 &
16597280453022 &
20644227272244012 
\\ \hline
$\frac{1}{125}h_7$ &
$-\frac{1}{5}$ & 13 & 2860 & 1855775 & 1750773750 & 1981335668498 & 2502724752660128 
\\ \hline
$\frac{1}{10}h_8$ &
$-\frac{1}{50}$ & 13 & 6425 & 6744325 & 8719953625 & 12525150549888 & 19171976431076873 
\\ \hline 
$\frac{1}{10}h_9$ &
$-\frac{1}{10}$ & 17 & 11185 & 12261425 & 16166719625 & 23478405649152 & 36191848368238417 \\ \hline
\end{tabular}

\end{center}
}

\section{Mirror quintic Calabi-Yau varieties}
\label{section1}
Let $W_{\psi}$ be the variety  obtained by a resolution of singularities of the following quotient:
\begin{equation}
 \label{shahva}
W_\psi:=\{[x_0:x_1:x_2:x_3:x_4]\in \mathbb P ^ 4\mid x_0^5+  x_1^5+ x_2^5+ x_3^5+ x_4^5-5\psi x_0x_1x_2x_3x_4=0\}/G, \ 
\end{equation}
where $G$ is the group 
$$
G:=\{(\zeta_1,\zeta_2,\cdots,\zeta_5)\mid  \zeta_i^5=1, \ \zeta_1\zeta_2\zeta_3\zeta_4\zeta_5=1 \}
$$  
acting in a canonical way, see for instance \cite{can91} and  $\psi^5\not=1$. 
The family $W_\psi$ is Calabi-Yau and it is mirror dual to the universal family of quintic varieties in 
$\Pn 4$. 
From now on we denote by $W$ such a variety and we call it the mirror quintic Calabi-Yau threefold. 
This section is dedicated to algebraic-geometric aspects, such as
moduli space and algebraic de Rham cohomology, of mirror quintic Calabi-Yau threefolds. We will use the algebraic de Rham cohomology $H^3_\dR(W)$ which is even defined 
for $W$ defined over an arbitrary field of characteristic zero. 
The original text of Grothendieck \cite{gro66} is still the main source of information for algebraic de Rham cohomology.
In the present text by the moduli of the objects $x$ we mean the  set of all $x$ quotiented by natural isomorphisms.

\subsection{Moduli space, I}
We first construct explicit affine coordinates for the moduli $S$ of the pairs $(W,\omega)$, where $W$ is a mirror quintic Calabi-Yau threefold and $\omega$ is a 
holomorphic differential 3-form on $W$. We have 
$$
S\cong \C^2\backslash\{(t_0^5-t_4)t_4=0\},
$$ 
where for $(t_0,t_4)$ we associate the pair $(W_{t_0,t_4},\omega_1)$. In the affine coordinates $(x_1,x_2,x_3,x_4)$, that is $x_0=1$,  $W_{t_0,t_4}$ is  given by
\begin{eqnarray*}
 W_{t_0,t_4} &:=& \{ f(x)=0\}/G,\\ 
f(x) &:=& -t_4-x_1^5-x_2^5-x_3^5-x_4^5+5t_0x_1x_2x_3x_4,
\end{eqnarray*}
and 
$$
\omega_1:= \frac{ dx_1\wedge dx_2\wedge dx_3\wedge dx_4}{df}.
$$
The multiplicative group $G_m:=\C^*$ acts on $S$ by:
 $$
(W,\omega)\bullet k=(W,k^{-1}\omega),\ k\in G_m,\ (W,\omega)\in S.
$$
In coordinates $(t_0,t_4)$ this corresponds to:
\begin{equation}
\label{poloar}
(t_0,t_4)\bullet k=(kt_0,k^{5}t_4),\ (t_0,t_4)\in S,\ k\in G_m.
\end{equation}
Two Calabi-Yau varieties in the family (\ref{shahva}) are isomorphic if and only if they have the same $\psi^5$. This and
(\ref{poloar}) imply that distinct pairs $(t_0,t_4)$ give non-isomorphic pairs $(W,\omega)$.  

\subsection{Gauss-Manin connection}
\label{gmI}
For a proper smooth family $W/T$ of algebraic varieties defined over a field $\k$ of characteristic zero, we have
the Gauss-Manin connection
$$
 \nabla:H_{\dR}^{i}(W/T)\to \Omega_T^1\otimes_{\O_T}H_{\dR}^{i}(W/T),
$$
where $H_{\dR}^{i}(W/T)$ is the $i$-th relative de Rham cohomology and 
$\Omega_T^1$ is the set of differential 1-forms on $T$. For simplicity we have assumed that $T$ is affine and $H_{\dR}^{i}(W/T)$ is a $\O_T$-module,
where $\O_T$ is the $\k$-algebra of regular function on $T$.
By definition of
a connection, $\nabla$ is $\k$-linear and satisfies the Leibniz rule
$$
\nabla(r\omega)=dr\otimes \omega+r\nabla\omega, \omega \in H_{\dR}^{i}(W/T),\ r\in \O_T.
$$
For a vector field $v$ in $T$ we define
$$
\nabla_v:  H_\dR^i(X)\to H_\dR^i(X) 
$$
to be $\nabla$ composed with 
$$
v\otimes {\rm Id}: \Omega^1_T\otimes_{\O_T} H_\dR^i(X)\to
\O_T\otimes_{\O_T} H_\dR^i(X)=H_\dR^i(X).
$$
Sometimes it is useful to choose   a basis $\omega_1,\omega_2,\ldots,\omega_h$ of the $\O_T$-module $H^i(X/T)$ and 
write the Gauss-Manin connection in this basis:
\begin{equation}
\label{18oct2010}
\nabla\begin{pmatrix}\omega_1 \\ \omega_2 \\ \vdots \\ \omega_h\end{pmatrix}=A\otimes \begin{pmatrix}\omega_1 \\ \omega_2 \\ \vdots \\ \omega_h\end{pmatrix}
\end{equation}
where $A$ is a $h\times h$ matrix with entries in $\Omega^1_T$. For further information on Gauss-Manin connection see \cite{kaod68}. See also \cite{ho06-1} for computational aspects
of Gauss-Manin connection.

\subsection{Intersection form and Hodge filtration}
\label{intersection}
For $\omega,\alpha\in H^3_\dR(W_{t_0,t_{4}})$ let
$$
\langle \omega,\alpha \rangle:={\rm Tr}(\omega\cup \alpha) 
$$
be the intersection form. 
If we consider $W$ as a  complex manifold and its de Rham cohomology defined by $C^\infty$ forms, then the intersection
form is just $\langle \omega,\alpha \rangle=\frac{1}{(2\pi i)^3}\int_{W_{t_0,t_{4}}}\omega\wedge \alpha$.
Using Poincar\'e duality it can be seen that it is dual  to the topological intersection form in $H_3(W_{t_0,t_4},\Q)$, for all these see for instance 
Deligne's lectures in \cite{dmos}. 
In  $H^{3}_\dR(W_{t_0,t_4})$  we have the Hodge filtration
$$
\{0\}=F^4\subset F^3\subset F^2\subset F^1\subset F^0=H^{3}_\dR(W_{t_0,t_4}),\ \ \dim_\C(F^i)=4-i.
$$
There is a relation between the Hodge filtration and the intersection form which is given by the following collection of equalities:
$$
\langle F^i,F^j\rangle=0, \ i+j\geq 4. 
$$
The Griffiths transversality is a property combining the Gauss-Manin connection and the Hodge filtration.
It says that  the Gauss-Manin connection sends $F^i$ to $\Omega^1_{S}\otimes F^{i-1}$ for $i=1,2,3$. 
Using this we conclude 
that:
\begin{equation}
\label{29oct11}
\omega_i:= {\frac{\partial^{i-1}}{\partial t_0^{i-1}}}(\omega_1)\in F^{4-i}, \ i=1,2,3,4.
\end{equation}
By abuse of notation we have used $\frac{\partial }{\partial t_0}$ 
instead of $\nabla_{\frac{\partial}{\partial t_0}}$.
 The intersection form in the basis $\omega_i$ is:

$$
 [\langle \omega_i,\omega_j\rangle]=
\begin{pmatrix}
0             &       0     &      0    &       \frac{1}{625}(t_4-t_0^5)^{-1} \\
0             &       0     &    -\frac{1}{625}(t_4-t_0^5)^{-1}  & -\frac{1}{125}t_0^4(t_4-t_0^5)^{-2}  \\         
0             &       \frac{1}{625}(t_4-t_0^5)^{-1}&0&           \frac{1}{125}t_0^3(t_4-t_0^5)^{-2} \\      
-\frac{1}{625}(t_4-t_0^5)^{-1} & \frac{1}{125}t_0^4(t_4-t_0^5)^{-2} & -\frac{1}{125}t_0^3(t_4-t_0^5)^{-2} &0      
\end{pmatrix}.
$$
For a proof see \cite{ho21}, page 468.

\subsection{Moduli space, II}
\label{modulispaceii}
We make the base change $\alpha=S\omega$, where $S$ is given by
\begin{equation}
\label{26apr2011}
S=\begin{pmatrix}
1 & 0& 0 & 0\\
t_{9} &
t_8 & 0 &0 \\
t_7
&
t_6
&
t_5
&
0 \\
t_1&t_2&t_3& t_{10}
\end{pmatrix}
\end{equation}
and $t_i$'s are unknown parameters, and we assume the the intersection form in $\alpha_i$'s is given by the matrix $\Phi$ in (\ref{31aug10}):
$$
\Phi=[\langle \alpha_i,\alpha_j\rangle]=S[\langle \omega_i,\omega_j\rangle ]S^\tr.
$$
This yields to many polynomial relations between $t_i$'s. It turns out that we can take $t_i,\ i=1,2,3,\ldots,6$ as
independent parameters and calculate all others in terms of these seven parameters:
\begin{eqnarray*}
t_7t_8-t_6t_9 &=& 3125t_0^3+t_2, \\ 
t_{10}&=& -t_8t_5,\\
t_5t_{9} &=& -3125t_0^4-t_3,\\
t_{10} &=& 625(t_4-t_0^5). \\
\end{eqnarray*}
The expression of $t_7,t_8,t_9$ are used in Theorem \ref{main2}.
For the moduli space $T$ introduced in Theorem \ref{main3} we get an isomorphism of sets (\ref{thanksdeligne}), 
where for $t$ in the right hand side of the  isomorphism (\ref{thanksdeligne}), we associate the pair $(W_{t_0,t_4},\alpha)$.
We also define
$$
\tilde t_5=\frac{1}{3125}t_5, \  \  \ \tilde t_6= -\frac{1}{5^6}(t_0^5-t_4)t_6+\frac{1}{5^{10}}(9375t_0^4+2t_3)t_5
$$
which correspond to the parameters $t_5,t_6$ in the previous article \cite{ho21}. 


\subsection{The Picard-Fuchs equation}
Let us consider the one parameter family of Calabi-Yau varieties $W_{t_0,t_4}$ with $t_0=1$ and $t_{4}=z$ and denote by $\eta$ the restriction of $\omega_1$ to these parameters. 
We calculate  the Picard-Fuchs equation of $\eta$ with respect to the parameter $z$:
$$
{\frac{\partial^{4}\eta}{\partial z^{4}}}=\sum_{i=1}^{4} a_i(z){\frac{\partial^{i-1}\eta}{\partial z^{i-1}}} \ \ \ \ \text{    modulo relatively exact forms.}
$$
This is in fact the linear differential equation 
\begin{equation}
\label{18fev2009}
I''''=\frac{-24}{625z^3(z-1)}I+
\frac{-24z+5}{5z^3(z-1)}I'+
\frac{-72z+35}{5z^2(z-1)}I''+
\frac{-8z+6}{z(z-1)}I''',\ \ '=\frac{\partial}{\partial z}
\end{equation}
which is calculated in \cite{can91}. This differential equation can be also written in the form (\ref{4jul2011}).
In \cite{ho06-1} we have developed algorithms which calculate such differential equations. The parameter $z$ is more convenient
for our calculations than the parameter $\psi$ and this is the reason why in this section we have used $z$ instead of $\psi$.  
The differential equation (\ref{18fev2009}) is satisfied by the periods 
$$
I(z)=\int_{\delta_z}\eta,\ \delta\in H_3(W_{1,z},\Q)
$$ 
of the differential form $\eta$ on the 
family $W_{1,z}$. 
In the basis ${\frac{\partial^{i}\eta}{\partial z^{i}}},\ \ i=0,1,2,3$ of $H^3_\dR(W_{1,z})$ the Gauss-Manin 
connection matrix has the form
\begin{equation}
\label{06jul2010}
A(z)dz:=
\begin{pmatrix}
0&1&0&0\\
0&0&1&0\\
0&0&0&1\\
a_1(z)&a_{2}(z)&a_{3}(z)&a_{4}(z)\\
\end{pmatrix}dz.
\end{equation}

\subsection{Gauss-Manin connection I}
We would like to calculate the Gauss-Manin connection 
$$
 \nabla:H_{\dR}^{3}(W/S)\to \Omega_S^1\otimes_{\O_S}H_{\dR}^{3}(W/S)
$$
of the  two parameter proper  family of 
varieties $W_{t_0,t_4},\ (t_0,t_4)\in S$.
 We calculate $\nabla$ with respect to the 
basis (\ref{29oct11}) 
 of $H^3_\dR(W/S)$. 
For this purpose we return back to the one parameter case. 
Now, consider the identity map 
$$
g:W_{{(t_0,t_{4})}}\to W_{1,z},\  
$$
which satisfies $g^*\eta=t_0\omega_1$. Under this map
$$
\frac{\partial}{\partial z}=\frac{-1}{5}\frac{t_0^{6}}{t_{4}}\frac{\partial}{\partial t_0}\left (=t_0^{5}\frac{\partial}{\partial t_{4}}\right ).
$$
From these two equalities we obtain a matrix $\tilde S=\tilde S(t_0,t_{4})$ such that
$$
[\eta, {\frac{\partial\eta}{\partial z}}, {\frac{\partial^2\eta}{\partial z^2}}, 
{\frac{\partial^{3}\eta_1}{\partial z^{3}}}]^\tr=
\tilde S^{-1}[\omega_1,\omega_2,\omega_3,\omega_{4}]^\tr, 
$$
where $\tr$ denotes the transpose of matrices, and the Gauss-Manin connection in the basis $\omega_i,\ i=1,2,3,4$ is:
$$
\tilde \GM=\left (d\tilde S+\tilde S\cdot A(\frac{t_{4}}{t_0^{5}})\cdot d(\frac{t_{4}}{t_0^{5}})\right )\cdot \tilde S^{-1},
$$
which is the following matrix after doing explicit calculations:
{\tiny
\begin{equation}
 \label{25aug2010}
\begin{pmatrix}
-\frac{1}{5t_4}dt_4
& dt_0+\frac{-t_0}{5t_4}dt_4
&0
&0
\\0
&\frac{-2}{5t_4}dt_4
&dt_0+\frac{-t_0}{5t_4}dt_4
&0
\\0
&0
&\frac{-3}{5t_4}dt_4
&dt_0+\frac{-t_0}{5t_4}dt_4
\\\frac{-t_0}{t_0^5-t_4}dt_0+\frac{t_0^2}{5t_0^5t_4-5t_4^2}dt_4
&\frac{-15t_0^2}{t_0^5-t_4}dt_0+\frac{3t_0^3}{t_0^5t_4-t_4^2}dt_4
&\frac{-25t_0^3}{t_0^5-t_4}dt_0+\frac{5t_0^4}{t_0^5t_4-t_4^2}dt_4
&\frac{-10t_0^4}{t_0^5-t_4}dt_0+\frac{6t_0^5+4t_4}{5t_0^5t_4-5t_4^2}dt_4
\end{pmatrix}.
\end{equation}
}
Now, we calculate the Gauss-Manin connection matrix of the family $W/T$ written in the basis $\alpha_i,\ i=1,2,3,4$. This is
$$
\GM=\left (dS+S\cdot \tilde {\GM}\right )\cdot S^{-1},
$$
where $S$ is the base change matrix (\ref{26apr2011}). Since the matrix $\GM$ is huge and does not fit into a mathematical paper, we do not
write it here.

\subsection{Modular differential equation}
We are in the final step of the proof of Theorem \ref{main3}. We have calculated the Gauss-Manin connection $\GM$ written in the basis 
$\alpha_i,\ i=1,2,3,4$.  It is a matter of explicit linear algebra  calculations to show that there is a unique vector field 
$\Ra$ in $T$ with the properties mentioned in Theorem (\ref{main3}) and to calculate it. 
In summary, the  Gauss-Manin connection matrix composed with the vector field $\Ra$ and written in the basis $\alpha_i$ has the form: 
\begin{equation}
 \label{17apr}
\nabla_{{\Ra}}=
\begin{pmatrix}
 0&1&0&0\\
0&0& \frac{5^8(t_4-t_0^5)^2}{t_5^3}&0\\
0&0&0&-1\\
0&0&0&0
\end{pmatrix}.
\end{equation}
It is interesting that the Yukawa coupling appears as the only non constant term in the above matrix.

\subsection{Algebraic group}
There is an algebraic group which acts on the right hand side of the isomorphism (\ref{thanksdeligne}). 
It corresponds to the base change in $\alpha_i,\ i=1,2,3,4$ such that the new
basis is still compatible with the Hodge filtration and we have still the intersection matrix (\ref{31aug10}):
$$
G:=\{ g=[g_{ij}]_{4\times 4}\in \GL (4,\C)\mid g_{ij}=0, \ \hbox{ for } j<i \hbox { and }g^\tr\Phi g=\Phi \},
$$
$$
\left \{g=\begin{pmatrix}
 g_{11}&g_{12}&g_{13}&g_{14}\\
0&g_{22}&g_{23}&g_{24}\\
0&0&g_{33}&g_{34}\\
0&0&0&g_{44}
\end{pmatrix}, \ g_{ij}\in\C    
\begin{matrix} 
g_{11}g_{44}=1,\\ 
g_{22}g_{33}=1,\\ 
g_{12}g_{44}+g_{22}g_{34}=0,\\ 
g_{13}g_{44}+g_{23}g_{34}-g_{24}g_{33}=0,\\
\end{matrix} 
\right \}.
$$
$G$ is called the Borel subgroup of ${\rm Sp}(4,\C)$ respecting the Hodge flag.
The action of $G$ on the moduli $T$ is given by:
$$
(W,[\alpha_1,\alpha_2,\alpha_3,\alpha_4])\bullet g= (W,[\alpha_1,\alpha_2,\alpha_3,\alpha_4]g).
$$
The algebraic group $G$ is of dimension six and  has two multiplicative subgroup $G_m=(\C^*,\cdot)$ and four additive subgroup $G_a=(\C,+)$ which generate it. 
In fact, an element $g\in G$ can be written in a unique way as the following product:
$$
\begin{pmatrix}
g_1^{-1}&-g_3g_1^{-1}&(-g_3g_6+g_4)g_1^{-1}&(-g_3g_4+g_5)g_1^{-1}\\
0&      g_2^{-1}     &g_6g_2^{-1}         &g_4g_2^{-1}\\         
0&      0           &g_2               &g_2g_3\\           
0&      0           &0                   &g_1  
 \end{pmatrix}=
$$
{\tiny
$$
\begin{pmatrix}
 g_1^{-1} &0&0&0\\
0&1&0&0\\
0&0&1&0\\
0&0&0&g_1
\end{pmatrix}
\begin{pmatrix}
 1 &0&0&0\\
0&g_2^{-1}&0&0\\
0&0&g_2&0\\
0&0&0&1
\end{pmatrix}
\begin{pmatrix}
 1 &-g_3&0&0\\
0&1&0&0\\
0&0&1&g_3\\
0&0&0&1
\end{pmatrix}
\begin{pmatrix}
 1&0&g_4&0\\
0&1&0&g_4\\
0&0&1&0\\
0&0&0&1
\end{pmatrix}
\begin{pmatrix}
 1 &0&0&g_5\\
0&1&0&0\\
0&0&1&0\\
0&0&0&1
\end{pmatrix}
\begin{pmatrix}
 1 &0&0&0\\
0&1&g_6&0\\
0&0&1&0\\
0&0&0&1
\end{pmatrix}.
$$
}
In other words, we have a bijection of sets $G_{m}\times G_{m}\times G_a\times G_a\times G_a\times G_a\cong G$ sending $(g_i)_{i=1,\ldots,6}$ to the above product.
If we identify an element $g\in G$ with the vector $(g_i)_{i=1,\ldots,6}$ then
$$
(g_1,g_2,g_3,g_4,g_5,g_6)^{-1}=
$$
$$
(g_1^{-1},g_2^{-1},-g_1^{-1}g_2g_3,  g_1^{-1}g_{2}^{-1}(g_3g_6-g_4), g_1^{-2}(-g_3^2g_6+2g_3g_4-g_5), -g_2^{-2}g_6 ).
$$
We denote by $\bullet$ the right action of  $G$ on the $t=(t_0,t_1,\ldots,t_6)$ space.
\begin{prop}
\label{19july2011}
The action of $G$ on $t_i$ (as a regular function on the affine variety $T$) is given by:
\begin{eqnarray*}
g\bullet t_0 &=& t_0g_1, \\
g\bullet t_1 &=&   t_1g_1^2+t_7g_1g_2g_3+t_{9}g_1g_2^{-1}g_4-g_3g_4+g_5,\\
g\bullet t_2  &=&  t_2g_1^3+t_6g_1^2g_2g_3+t_8g_1^2g_2^{-1}g_4,\\
 g\bullet t_3  &=& t_3g_1^4+t_5g_1^3g_2g_3,\\
g\bullet t_4  &=& t_4g_1^5,\\
g\bullet t_5 &=&   t_5g_1^3g_2,\\
g\bullet t_6 &=&  t_6g_1^2g_2+t_8g_1^2g_2^{-1}g_6.\\
\end{eqnarray*}
Consequently 
\begin{eqnarray*}
g\bullet t_7 &=&   t_7g_1g_2+t_{9}g_1g_2^{-1}g_6-g_3g_6+g_4,\\
g\bullet t_8 &=&      t_8g_1^2g_2^{-1},\\
g\bullet t_{9} &=&   t_{9}g_1g_2^{-1}-g_3, \\
g\bullet t_{10} &=&  t_{10}g_1^5.
\end{eqnarray*}

\end{prop}


\begin{proof}
We first calculate the action of $g=(k,1,0,0,0,0),\ k\in \C^*$ on $t$.
We have an isomorphism $(W_{(t_0,t_4)},k^{-1}\omega_1)\cong (W_{(t_0k,t_4k^{5})},\omega_1)$ given by 
$$
(x_1,x_2,x_3,x_4)\mapsto (k^{-1}x_1,k^{-1}x_2,k^{-1}x_3,k^{-1}x_4).
$$
Under this isomorphism the vector field $k^{-1}\frac{\partial}{\partial t_0}$ is mapped to $\frac{\partial}{\partial t_0}$ and so $k^ {-i}\omega_i$ is mapped to $\omega_i$.
This implies the isomorphisms
$$
(W_{(t_0,t_4)},k(t_1\omega_1+t_2\omega_2+t_3\omega_3+625(t_4-t_0^5)\omega_4))\cong 
$$
$$(W_{(t_0,t_4)\bullet k},k^2t_1\omega_1+k^3t_2\omega_2+k^4t_3\omega_3+625(k^5t_4-(kt_0)^5)\omega_4))
$$
and
$$
(W_{(t_0,t_4)},S\omega )\cong (W_{(t_0,t_4)\bullet k},S\begin{pmatrix}
 k &&0&0\\
0&k^2&0&0\\
0&0&k^3&0\\
0&0&0&k^4
\end{pmatrix}  \omega),
$$
where $S$ is defined in (\ref{26apr2011}). Therefore, 
$$
g\bullet t_i=t_ik^{\tilde d_i}, \ \tilde d_i=i+1, \ i=0,1,2,3,4\ \ \ \tilde d_5=3, \ \ \tilde d_6=2  
$$

\end{proof}

\section{Periods}
\label{section2}
This section is dedicated to transcendental aspects of mirror quintic Calabi-Yau threefold. By this we mean the periods
of meromorphic differential 3-forms over topological cycles. We first work with periods without calculating them explicitly.

\subsection{Period map}
\label{1111}
We choose a symplectic basis for $H_3(W,\Z)$, that is, a basis $\delta_i,\ i=1,2,3,4$ such that
$$
\Psi:=[\langle \delta_i,\delta_j \rangle]=\begin{pmatrix}
 0&0&1&0\\
0&0&0&1\\
-1&0&0&0\\
0&-1&0&0
\end{pmatrix}.
$$
It is also convenient  to use the basis 
$$
[\tilde \delta_1,\tilde \delta_2,\tilde \delta_3, \tilde \delta_4]=[\delta_1,\delta_2,\delta_3, \delta_4]\Psi^{-1}=[\delta_3,\delta_4,-\delta_1,-\delta_2].
$$ 
In this basis the intersection 
form is $[\langle \tilde \delta_i,\tilde \delta_j\rangle]=\Psi^{-\tr}=\Psi$. 
Recall that in \S\ref{modulispaceii} a mirror quintic Calabi-Yau threefold  $W$ is equipped with a basis $\alpha_1,\alpha_2,\alpha_3,\alpha_4$ of $H_\dR^3(W)$ compatible with the Hodge 
filtration  and such that 
$[\langle \alpha_i,\alpha_j\rangle]=\Phi$. We define the period matrix to be
$$
[x_{ij}]=[\int_{\delta_i} \alpha_j].
$$
(In this section we discard the usage of $x_{ij}$ in the Introduction). 
Let $\tilde \delta_i^p\in H^3(W,\Q)$ be the Poincar\'e dual of $\tilde \delta_i$, that is, 
it is defined by the property 
$\int_{\delta}\tilde \delta_i^p=\langle \tilde\delta_i,\delta\rangle$ for all $\delta\in H_3(W,\Q)$. If we write 
$\alpha_i$ in terms of $\tilde \delta_i^p$ what we get is:
$$
[\alpha_1,\alpha_2,\alpha_3,\alpha_4]=[\tilde \delta_1^p,\tilde \delta_2^p,\tilde \delta_3^p, \tilde \delta_4^p][\int_{\delta_i}\alpha_j],
$$    
that is, the coefficients of the base change matrix are the periods of $\alpha_i$'s over $\delta_i$'s and not $\tilde \delta_i$'s.
We have 
\begin{equation}
 \label{24aug10}
[\langle \alpha_i,\alpha_j\rangle]=[\int_{\delta_i}\alpha_j]^{\tr} \Psi^ {-\tr}[\int_{\delta_i}\alpha_j] .
\end{equation}
and so we get:
\begin{equation}
\label{6oct2011}
\Phi-[x_{ij}]^{\tr}\Psi [x_{ij}]=0.
\end{equation}
This gives us  6 non trivial polynomial relations between periods $x_{ij}$:
\begin{eqnarray*}
x_{12}x_{31}-x_{11}x_{32}+x_{22}x_{41}-x_{21}x_{42}\ \ \ \ \   &= &0,\\ 
x_{13}x_{31}-x_{11}x_{33}+x_{23}x_{41}-x_{21}x_{43}\ \ \ \ \  & =& 0,\\
x_{14}x_{31}-x_{11}x_{34}+x_{24}x_{41}-x_{21}x_{44}+1 &=& 0,\\ 
x_{13}x_{32}-x_{12}x_{33}+x_{23}x_{42}-x_{22}x_{43}+1 &=& 0,\\
x_{14}x_{32}-x_{12}x_{34}+x_{24}x_{42}-x_{22}x_{44}\ \ \ \ \  &=& 0,\\ 
x_{14}x_{33}-x_{13}x_{34}+x_{24}x_{43}-x_{23}x_{44}\ \ \ \ \  &=& 0.
\end{eqnarray*}
These equalities correspond to the entries $(1,2),(1,3),(1,4),(2,3),(2,4)$ and  $(3,4)$ of (\ref{6oct2011}). 
Taking the determinant of (\ref{6oct2011}) we see that up to sign we have $\det(\per)=-1$. 
There is another effective way to calculate this determinant without
the sign ambiguity. In the ideal of $\Q[x_{ij}, \ i,j=1,2,3,4]$
generated by the polynomials $f_{12},f_{13},f_{14},f_{23},f_{2,4},f_{34}$ in the right hand side of the above equalities, 
the polynomial $\det([x_{ij}])$ is reduced to $-1$. 
Let $y_{ij}$ be indeterminate variables, $R=\C[y_{ij}, i,j=1,2,3,4]$ and  $I:=\{f\in R\mid f(\cdots, x_{ij},\cdots)=0\}$.
The ideal $I$ is generated by $f_{12},f_{13},f_{14},f_{23},f_{2,4},f_{34}$, see for instance \cite{ho21}, Proposition 3 page 472.

\subsection{A special locus}
\label{sI}
Let 
$$
C^\tr:=[0,1,0,0][\langle \tilde \delta_i,\tilde\delta_j\rangle]^{-\tr}=[0,0,0,1].
$$ 
We are interested in the loci $L$ of parameters $t\in T$ such that
\begin{equation}
\label{badbad}
[\int_{\delta_1}\alpha_4,\ldots, \int_{\delta_4}\alpha_4]=C.
\end{equation}
Using the equality corresponding to the $(1,4)$ entries of (\ref{24aug10}), we note that on this locus 
we have
$$
\int_{\delta_2}\alpha_1=1,\ \ \int_{\delta_2}\alpha_i=0,\ i\geq 2.
$$
The equalities (\ref{badbad}) define a three dimensional locus of $T$. We also put the following two conditions
$$
\int_{\delta_1}\alpha_2=1,\ \int_{\delta_1}\alpha_3=0
$$
in order to get a dimension one locus. Finally using (\ref{6oct2011}) we conclude that the period matrix 
for points in $L$ is of the form
\begin{equation}
 \tau=\label{UERJ?}
\begin{pmatrix}
 \tau_0 &1&0&0\\
1&0&0&0\\
\tau_{1}& \tau_{3}& 1&0\\
\tau_{2} & -\tau_0 \tau_{3}+\tau_{1} & -\tau_0& 1 
\end{pmatrix}.
\end{equation}
The particular expressions for the $(4,2)$ and $(4,3)$ entries of the above matrix follow from the polynomial relations  (\ref{6oct2011}). 
The Gauss-Manin connection matrix restricted to $L$ is:
$$ 
\GM\mid_L=d\tau^\tr\cdot \tau^{-\tr}=
\begin{pmatrix}
0&d\tau_0& -\tau_{3}d\tau_0+d\tau_{1} &  -\tau_{1}d\tau_0+\tau_0 d\tau_{1}+d\tau_{2}\\  
0&0&d\tau_{3}&-\tau_{3}d\tau_0+d\tau_{1} \\ 
0&0&       0&       -d\tau_0\\
0&0&       0&       0 
\end{pmatrix}.        
$$
The Griffiths transversality theorem implies that 
$$
-\tau_{3}d\tau_0+d\tau_{1}=0,\  \ -\tau_{1}d\tau_0+\tau_0 d\tau_{1}+d\tau_{2}=0.
$$
Since $L$ is one dimensional, there are analytic relations between $\tau_i,\ i=0,1,2,3$. So, we consider $\tau_0$ as 
an independent parameter and 
$\tau_{1},\tau_{2},\tau_{3}$ depending on $\tau_0$, then we have 
\begin{equation}
\label{retrovisor}
\tau_{3}=\frac{\partial \tau_1}{\partial \tau_0},\ \frac{\partial\tau_2}{\partial \tau_0}=\tau_{1}-\tau_0 \frac{\partial\tau_1}{\partial \tau_0}. 
\end{equation}
We conclude that the Gauss-Manin connection matrix restricted to $L$ and composed with the vector field $\frac{\partial}{\partial \tau_0}$ is given by:
\begin{equation}
\label{17apr2011}
\begin{pmatrix}
0&1&0 &  0\\  
0&0&\frac{\partial \tau_{3}}{\partial \tau_0}&0 \\ 
0&0&       0&       -1\\
0&0&       0&       0 
\end{pmatrix}.        
\end{equation}
 
\begin{prop}
The functions  $t_i(\tau_0)$ obtained by the regular functions $t_i,\ i=0,1,2,\ldots,6$ restricted to $L$ and seen as functions in 
$\tau_0$ form a solution of the ordinary differential equation $\Ra$.  
\end{prop}
\begin{proof}
 It follows from (\ref{17apr2011})  and the uniqueness of the vector field $\Ra$ satisfying the equalities (\ref{17apr}).
\end{proof}




\subsection{The algebraic group and special locus}
\label{sarasara}
For any $4\times 4$ matrix $x=[x_{ij}]$ satisfying (\ref{6oct2011}) and  
\begin{equation}
\label{10oct2011}
x_{11}x_{22}-x_{12}x_{21}\not =0,\ x_{21}\not=0,  
\end{equation}
 there is a unique $g\in G$ such $xg$ is of the form (\ref{UERJ?}). To prove this affirmation 
explicitly, we take an arbitrary $x$ and $g$ and we write down the corresponding equations corresponding to the six 
entries $(2,1),(1,2), (2,2), (1,3),(2,3), (2,4)$ of $xg$, that is
$$
xg=
\begin{pmatrix}
 *&1& 0&* \\
1& 0& 0& 0\\
*&*&*&*\\
*&*&*&*
\end{pmatrix}.
$$
For our calculations we will need the coordinates of $g^{-1}$ in terms of $x_{ij}$:
\begin{eqnarray*}
g_1 &=& x_{21}^{-1},\\
g_2 &=& \frac{-x_{21}}{x_{11}x_{22}-x_{12}x_{21}},\\
g_3 &=& \frac{-x_{22}}{x_{21}},\\
g_4 &=& \frac{-x_{12}x_{23}+x_{13}x_{22}}{x_{11}x_{22}-x_{12}x_{21}},\\
g_5 &=& \frac{x_{11}x_{22}x_{24}-x_{12}x_{21}x_{24}+x_{12}x_{22}x_{23}-x_{13}x_{22}^2}{x_{11}x_{21}x_{22}-x_{12}x_{21}^2},\\
g_6 &=& \frac{x_{11}x_{23}-x_{13}x_{21}}{x_{11}x_{22}-x_{12}x_{21}}.
\end{eqnarray*}
Substituting the expression of $g$ in terms of $x_{ij}$ in $\tau=xg$ we get:
$$
\tau=\begin{pmatrix}
\frac{x_{11}}{x_{21}}&1&0&0\\
1&0&0&0\\
\frac{x_{31}}{x_{21}}&  \frac{-x_{21}x_{32}+x_{22}x_{31}}{x_{11}x_{22}-x_{12}x_{21}}  & 1 & 0\\
\frac{x_{41}}{x_{21}} &  \frac{-x_{21}x_{42}+x_{22}x_{41}}{x_{11}x_{22}-x_{12}x_{21}}   & -\frac{x_{11}}{x_{21}} & 1
 \end{pmatrix}.
$$
Note that for the entries $(1,4),\ (3,3)$ and $(4,3)$ of the above matrix 
we have used the polynomial relations (\ref{6oct2011}) between periods. 

\section{Monodromy covering}
\label{18july2011}
\label{section3}
In the previous section we described a solution of $\Ra$ locally. In this section we study further such a solution in a global context.
More precisely, we describe a meromorphic map  $t:\uhp \to T$ whose image is $L$ of the previous section, where $\uhp$ is 
the monodromy covering of $(\ref{4jul2011})$. 

\subsection{Monodromy covering}
Let $\tilde \uhp$ be the moduli of the pairs $(W, \delta)$, where $W$ is a mirror quintic Calabi-Yau threefold 
and $\delta=\{\delta_1,\delta_2,\delta_3,\delta_4\}$ is a basis  of $H_3(W,\Z)$ such that the intersection matrix
in this basis is $\Psi$, that is, $[\langle \delta_i,\delta_j\rangle]=\Psi$.
The set $\tilde \uhp$ has a canonical structure of a Riemann surface, not necessarily connected.
We denote by $\uhp$ the connected component of $\tilde\uhp$ which contain the  particular pair $(W_{1,z},\delta)$
such that the monodromies around 
$z=0$ and $z=1$ are respectively given by the matrices $M_0$ and $M_1$ in the Introduction. 
It is conjectured that in the monodromy group $\Gamma:=\langle M_0,M_1\rangle$ the only relation between $M_0$ and $M_1$ is 
$(M_0M_1)^5=I$. This is equivalent to say that $\uhp$ is biholomorphic to the
upper half plane. We do not need to assume this conjecture. 
 By definition, the monodromy group
$\Gamma$ acts on $\uhp$ by base change in $\delta$.
The bigger group ${\rm Sp}(4,\Z)$ acts also on $\tilde \uhp$ by base change
and all connected components of $\tilde \uhp$ are obtained by $ \uhp_\alpha:=\alpha(\uhp),\ \alpha\in {\rm Sp}(4,\Z)/\Gamma$:
$$
\tilde \uhp:=\cup_{\alpha\in  {\rm Sp}(4,\Z)/\Gamma} \tilde \uhp_\alpha.
$$
From now on by $w$ we denote a point $(W,\delta)$ of $\uhp$. We use the following meromorphic functions on $\uhp$:
$$
\tau_i:\uhp\to \C,\ i=0,1,2,\ 
$$
$$
\tau_0(w)=\frac{\int_{\delta_{1}}\alpha_1}{ \int_{\delta_{2}}\alpha_1},\ 
\tau_1(w)=\frac{\int_{\delta_{3}}\alpha_1}{ \int_{\delta_{2}}\alpha_1},\ 
\tau_2(w)=\frac{\int_{\delta_{4}}\alpha_1}{ \int_{\delta_{2}}\alpha_1},
$$
where $\alpha_1$ is a holomorphic differential form on $W$. They do not depend on the choice of $\alpha_1$. For simplicity, 
we have used the same notations $\tau_i$ as in \S\ref{section2}.

 There is a useful meromorphic function $z$ on $\uhp$ which is obtained by identifying $W$ with some $W_{1,z}$. It has a pole of order $5$ at 
elliptic points which are the pairs $(W,\delta)$ with $W=W_{\psi,1},\ \psi=0$.  In this way, we have a well-defined holomorphic function 
$$
\psi=z^{-\frac{1}{5}}:\uhp\to \C.
$$
The coordinate system $\tau_0$ is  adapted for calculations around the cusp $z=0$.
 Let $B$ be the set of points $w=(W,\delta)\in \uhp$ such that either $\tau_0$ has a pole at $w$ or it has a 
critical point  at $w$, that is,  $\frac{\partial \tau_0}{\partial z}(w)=0$. 
 We do not know whether $B$ is empty or not. Many functions that we are going
to study are meromorphic with poles at $B$. The set $B$ is characterized by this property that 
in its complement in  $\uhp$ 
the inequalities (\ref{10oct2011}) hold.

\subsection{A particular solution}
\label{21july2011}
For a point $w=(W,\delta)\in\uhp\backslash B$ there is a unique basis $\alpha$ of $H^3_\dR(W)$ such that  $(W,\alpha)$ is an element in the moduli space $T$ defined
in \S\ref{modulispaceii} and the period matrix $[\int_{\delta_i}\alpha_j]$ of the triple $(W,\delta,\alpha)$ is of the form (\ref{UERJ?}). 
This follows from the arguments in \S\ref{sarasara}. In this way we have well-defined meromorphic maps
$$
t: \uhp\to T
$$
and
$$
\tau: \uhp\to \Mat(4,\C)
$$
which are characterized by the uniqueness of the basis $\alpha$ and the equality:
$$
\tau(w)=[\int_{\delta_i}\alpha_j].
$$
If we parameterize $\uhp$ by the image of $\tau_0$ then  $t$ is the same map as in \S\ref{sI}. We conclude that the map
$t:\uhp\to T$ with the coordinate system $\tau_0$ on $\uhp$ is a solution of $\Ra$. The functions $t$ and $\tau$ are holomorphic outside
the poles and critical points of $\tau_0$ (this corresponds to points in which the inequalities (\ref{10oct2011}) occur).

\subsection{Action of the monodromy}
The monodromy group $\Gamma:=\langle M_0,M_1\rangle$ acts on $\uhp$ by base change. If we choose the local 
coordinate system $\tau_0$ on $\uhp$ then this action is given by:
$$
A(\tau_0)=\frac{a_{11}\tau_0 +a_{12}+ a_{13}\tau_1+ 
a_{14}\tau_2}{a_{21}\tau_0 +a_{22}+ a_{23}\tau_1+ a_{24}\tau_2},\ \ A=[a_{ij}]\in\Gamma.
$$ 

\begin{prop}
 For all $A\in \Gamma$ we have
$$
t(w)=t(A(w))\bullet g(A,w),
$$
where $g(A,w)\in G$ is defined using the equality
$$
A\cdot \tau(w)=\tau(A(w))\cdot g(A,w).
$$
\end{prop}
\begin{proof}
Let $w=(W,\delta)\in\uhp$ and  $t(w)=(W,\alpha)$. By definition we have
$$
[\int_{A(\delta)_i}\alpha_j]g(A,w)^{-1}=A\tau(w)g(A,w)^{-1}=\tau(A(w)).
$$
Therefore, $t(A(w))=(W, \alpha\cdot g(A,w)^{-1})=t(w)\bullet g(A,w)^{-1}$.
\end{proof}
 
If we choose the coordinate system $\tau_0$ on $\uhp$ and regard the parameters $t_i$'s and $\tau_i$'s as functions in 
$\tau_0$, then we have
$$
t(\tau_0)=t(A(\tau_0))\bullet g(A,\tau_0).
$$
These are the functional equations of $t_i(\tau_0)$'s mentioned in the  Introduction. For $A=M_0$ we have:
{\tiny
$$
\begin{pmatrix}
1&1& 0& 0\\
0&1& 0& 0\\
5&5& 1& 0\\
0&-5&-1&1 
\end{pmatrix}
\begin{pmatrix}
 \tau_0 &1&0&0\\
1&0&0&0\\
\tau_1& \tau_{3}& 1&0\\
\tau_{2} & -\tau_0 \tau_{3}+\tau_{1} & -\tau_0& 1 
\end{pmatrix}=
\begin{pmatrix}
 \tau_0+1 &1&0&0\\
1&0&0&0\\
\tau_1+5\tau+5 & \tau_{3}+5& 1&0\\
\tau_{2}-5-\tau_1 & -\tau_0 (\tau_{3}+1) +\tau_{1} & -\tau_0-1& 1 
\end{pmatrix}
$$
}
which is already of the format (\ref{UERJ?}). Note that
$$
-(\tau_0+1)( \tau_{3}+5)+\tau_1+5\tau_0+5=-\tau_0 (\tau_{3}+1) +\tau_{1}.
$$
Therefore, $M_0(\tau_0)=\tau_0+1$ and $g(M_0,\tau_0)$ is the identity matrix. 
The corresponding functional equation of $t_i$ simply says that $t_i$ 
is invariant under $\tau_0\mapsto \tau_0+1$:
$$
t_i(\tau_0)=t_i(\tau_0+1),\ i=0,1,\ldots,6.
$$ 
For $A=M_1$ we have 
$$
\begin{pmatrix}
1&0&0&0\\
0&1&0&1\\
0&0&1&0\\
0&0&0&1 
\end{pmatrix}
\begin{pmatrix}
 \tau_0 &1&0&0\\
1&0&0&0\\
\tau_1& \tau_{3}& 1&0\\
\tau_{2} & -\tau_0 \tau_{3}+\tau_{1} & -\tau_0& 1 
\end{pmatrix}=
\begin{pmatrix}
 \tau_0 &1&0&0\\
\tau_{2}+1 & -\tau_0 \tau_{3}+\tau_{1} & -\tau_0& 1 \\
\tau_1& \tau_{3}& 1&0\\
\tau_{2} & -\tau_0 \tau_{3}+\tau_{1} & -\tau_0& 1 
\end{pmatrix}=
$$
{\tiny
$$
\begin{pmatrix}
\frac{\tau_0}{\tau_2+1} &1&0&0\\
1&
0&
0&
0\\
\frac{\tau_1}{\tau_2+1}&
\frac{\tau_0\tau_1\tau_3-\tau_1^2+\tau_2\tau_3+\tau_3}{\tau_0^2\tau_3-\tau_0\tau_1+\tau_2+1}&
1&
0\\
\frac{\tau_2}{\tau_2+1}&
\frac{-\tau_0\tau_3+\tau_1}{\tau_0^2\tau_3-\tau_0\tau_1+\tau_2+1}&
\frac{-\tau_0}{\tau_2+1}&
1
\end{pmatrix}
\begin{pmatrix}
 (\tau_2+1)&
(-\tau_0\tau_3+\tau_1)&
(-\tau_0)&
1\\
0&
\frac{\tau_0^2\tau_3-\tau_0\tau_1+\tau_2+1}{\tau_2+1}&
\frac{\tau_0^2}{\tau_2+1}&
\frac{-\tau_0}{\tau_2+1}\\
0&
0&
\frac{\tau_2+1}{\tau_0^2\tau_3-\tau_0\tau_1+\tau_2+1}&
\frac{\tau_0\tau_3-\tau_1}{\tau_0^2\tau_3-\tau_0\tau_1+\tau_2+1}\\
0&
0&
0&
\frac{1}{\tau_2+1}
\end{pmatrix},
$$
}
where the element of the algebraic group $G$ in the right hand side has the coordinates:
\begin{eqnarray*}
g_1 &=& \frac{1}{\tau_2+1},\\
g_2 &=& \frac{\tau_2+1}{\tau_0^2\tau_3-\tau_0\tau_1+\tau_2+1},\\
g_3 &=& \frac{\tau_0\tau_3-\tau_1}{\tau_2+1}, \\
g_4 &=& \frac{-\tau_0}{\tau_0^2\tau_3-\tau_0\tau_1+\tau_2+1},\\
g_5 &=& \frac{1}{\tau_0^2\tau_3-\tau_0\tau_1+\tau_2+1},\\
g_6 &=& \frac{\tau_0^2}{\tau_0^2\tau_3-\tau_0\tau_1+\tau_2+1}.
\end{eqnarray*}
In this case we have
$$
M_1(\tau_0)=\frac{\tau_0}{\tau_2+1}.
$$
The corresponding functional equations of $t_i$'s can be written immediately. These are presented in Theorem \ref{main2}.

%
%


\subsection{The solution in terms of periods}
In this section we explicitly calculate the map $t$.
For $w=(W,\delta)\in \uhp$ we identify $W$ with $W_{1,z}$ and hence we obtain a unique point 
$\tilde z=(1,0,0,0,z,1,0) \in T$. Now, 
we have a well-defined period map
$$
\per: \uhp\to \Mat(4,\C),
$$
$$
w=(W_{1,z}, \{\delta_{i,z},\ i=1,2,3,4\})\mapsto[\int_{\delta_{i,z}}\alpha_j].
$$
We write $\per(w)g(w)=\tau(w)$, where $\tau(w)$ is of the form (\ref{UERJ?}) and $g(w)\in G$. We have
$$
t(w)=\tilde z\bullet g(w).
$$
For the one dimensional locus $\tilde z\in T$, we have $\alpha=S\omega$ and $\omega=T\tilde \eta$, where
$$
S=\begin{pmatrix}
1&0&0&0\\
-5^5&-5^4(z-1)&0&0\\
-\frac{5}{z-1}&0&1&0\\
0&0&0&5^4(z-1)
\end{pmatrix}, \ \ 
T=\left(
\begin{array}{*{4}{c}}
1 & 0 & 0 & 0 \\
-1 & -5 & 0 & 0 \\
2 & 15 & 25 & 0 \\
-6 & -55 & -150 & -125
\end{array}
\right)
$$
and
$$
\tilde \eta=[\eta, \theta\eta,\theta^2\eta, \theta^3\eta ]^\tr,\   \ \theta=z\frac{\partial}{\partial z}.
$$
Therefore, $\alpha=ST\eta$.
Restricted to $\tilde z$-locus we have $\alpha_1=\omega_1=\eta$ and by our definition of $x_{ij}$'s in the introduction 
$$
x_{ij}=\theta^{j-1}\int_{\delta_i}\eta,\ i,j=1,2,3,4.
$$
Therefore,
$$
\per(w)=[x_{ij}](ST)^\tr.
$$
Now, the map $w\mapsto t(w)$, where the domain $\uhp$ is equipped with the coordinate system $z$, is given by 
the expressions for $t_i$ in Theorem \ref{main}. We conclude that
if we write $t_i$'s in terms of $\tau_0$  then we get functions which are solutions to
$\Ra$.  Note that
$$
\frac{\partial}{\partial \tau_0}=2\pi iq\frac{\partial }{\partial q}=
(z\frac{\partial\frac{x_{11}}{x_{21}}}{\partial z})^{-1}z\frac{\partial}{\partial z}=
\frac{x_{21}^2}{x_{12}x_{21}-x_{11}x_{22}}\theta.
$$

\subsection{Calculating periods}
\label{periodcalculation}
\label{section4}
In this section we calculate the periods $x_{ij}$ explicitly. 
This will finish the proof of our main theorems announced in 
the Introduction. 

We restrict the parameter $t\in T$ to the one dimensional loci $\tilde z$ given by $t_0=1, t_1=t_2=t_3=0, t_4=z,t_5=1,t_6=0$. 
On this locus $\eta=\omega_1=\alpha_1$. 
We know that the integrals $\int_\delta \eta,\ \delta\in H_3(W_{1,z},\Q)$ satisfy the linear differential equation (\ref{18fev2009}). 
Four linearly independent solutions of (\ref{18fev2009}) are given by 
$\psi_0,\psi_1,\psi_2,\psi_3$ in the Introduction, 
see for instance \cite{vanvan} and \cite{ChYaYu}.  
In fact, there are  four topological cycles with complex coefficients 
$\hat\delta_1,\hat\delta_2,\hat \delta_3,\hat\delta_4\in 
H_3(W_{1, z},\C)$ such that 
$\int_{\hat \delta_{i}} \eta=\frac{(2\pi i)^{i-1}}{5^4}\psi_{4-i}$.
Note that the pair $(W_{1,z},5\eta)$ is  isomorphic to the pair $(W_\psi,\Omega)$ used in \cite{can91}. 
 We use a new basis given by
$$
\begin{pmatrix}
 \delta_1\\ \delta_2 \\ \delta_3 \\ \delta_4
\end{pmatrix}=\begin{pmatrix}
   0 &0&1&0\\
0&0&0&1\\
0&d&\frac{d}{2}&-b\\
-d&0&-b&-a
  \end{pmatrix}
\begin{pmatrix}
\hat \delta_1\\ \hat \delta_2 \\ \hat\delta_3 \\ \hat\delta_4
\end{pmatrix},
$$
where
$$
a=\frac{c_3}{(2\pi i)^3}\zeta(3)=\frac{-200}{(2\pi i)^3}\zeta(3),\ \ b=c_2\cdot H/24=\frac{25}{12},\  d=H^3=5,
$$
(these notations are used in \cite{vanvan}). The monodromies around $z=0$ and $z=1$ written in the basis $\delta_i$ are respectively 
$$
M_0:=
\begin{pmatrix}
1&1& 0& 0\\
0&1& 0& 0\\
d&d& 1& 0\\
0&-k&-1&1 
\end{pmatrix}\ \ \ 
M_1:=
\begin{pmatrix}
1&0&0&0\\
0&1&0&1\\
0&0&1&0\\
0&0&0&1 
\end{pmatrix},
$$
where $k=2b+\frac{d}{6}=5$, see \cite{ChYaYu}. In fact, $\delta_i\in H_3(W_{1,z},\Z),\ i=1,2,3,4$. This follows from the  the calculations in 
\cite{can91} and the expressions for monodromy matrices. 
In summary, we have
\begin{eqnarray*}
\label{tanhayi1}
x_{11}=\int_{\delta_{1}}\eta&=&\frac{1}{5^2}(\frac{2\pi i}{5})^2\psi_1(\tilde z),\\
x_{21}=\int_{\delta_{2}}\eta&=& \frac{1}{5}(\frac{2\pi i}{5})^3\psi_0,\\
x_{31}=\int_{\delta_{3}}\eta&=& \frac{d}{125}\psi_2(\tilde z)\frac{2\pi i}{5}+\frac{d}{50}\cdot(\frac{2\pi i}{5})^2 \cdot\psi_1(\tilde z)-\frac{b}{5} \cdot(\frac{2\pi i}{5})^3 \cdot\psi_0(\tilde z),\\
x_{41}=\int_{\delta_{4}}\eta&=&\frac{-d}{5^4}\psi_3(\tilde z) +\frac{-b}{5^2}\cdot (\frac{2\pi i}{5})^2\cdot \psi_1(\tilde z) -\frac{a}{5}\cdot (\frac{2\pi i}{5})^3\cdot \psi_0(\tilde z), 
\end{eqnarray*}
where $\tilde z=\frac{z}{5^5}$. 
We have also
\begin{eqnarray*}
\tau_0=\frac{\int_{\delta_{1}}\eta}{\int_{\delta_{2}}\eta }&=&\frac{1}{2\pi i}\frac{\psi_1(\tilde z)}{\psi_0(\tilde z)},\\
\tau_1=\frac{\int_{\delta_{3}}\eta}{\int_{\delta_{2}}\eta } &=& d(\frac{1}{2} \tau_0^2+\frac{1}{5}H')+\frac{d}{2} \tau_0-b=-b+\frac{d}{2}\tau_0(\tau_0+1)+\frac{d}{5}H',\\
\tau_2=\frac{\int_{\delta_{4}}\eta}{\int_{\delta_{2}}\eta }&=&-d\left (\frac{-1}{3}\tau_0^3+\tau_0(\frac{1}{2}\tau_0^2+\frac{1}{5}H')+\frac{2}{5}H\right ) -b\tau_0 -a\\
&=& -a-b\tau_0-\frac{d}{6}\tau_0^3-\frac{d}{5}\tau_0H'-\frac{2d}{5}H,
\end{eqnarray*}
where $H$ is defined in (\ref{dastgir}). We have used the equalities
$$
\frac{\psi_2}{\psi_0}-\frac{1}{2}(\frac{\psi_1}{\psi_0})^2=\frac{1}{5}( \sum_{n=1}^\infty (\sum_{d|n}n_d d^3)\frac{q^n}{n^2}),
$$
$$
 \frac{1}{3}(\frac{\psi_1}{\psi_0})^3-\frac{\psi_1}{\psi_0}\frac{\psi_2}{\psi_0}+\frac{\psi_3}{\psi_0}=
\frac{2}{5} \sum_{n=1}^\infty (\sum_{d|n}n_d d^3)\frac{q^n}{n^3}
$$
see for instance \cite{kon95, pan97}. We can use the explicit series
\begin{eqnarray*}
\psi_0(\tilde z)&=&\sum_{m=0}^{\infty}\frac{(5m)!}{(m!)^5}\tilde z^m \\
\psi_1(\tilde z)&=& \ln(\tilde z)\psi_0(\tilde z)+5\tilde \psi_1(\tilde z),\ \ \tilde \psi_1(\tilde z)=\sum_{m=1}^{\infty}\frac{(5m)!}{(m!)^5}(\sum_{k=m+1}^{5m}\frac{1}{k})\tilde z^m 
\end{eqnarray*}
and calculate the $q$-expansion of $t_i(\tau_0)$ around the cusp $z=0$. 
There is another way of doing this using the differential equation $\Ra$. We just
use the above equalities to obtain the initial values (\ref{22july2010}) in the Introduction. 
We write each $h_i$ as a formal power series in $q$, $h_i=\sum_{n=0}^\infty t_{i,n}q^n$, and substitute in (\ref{lovely}) with $\dot t:=5q\frac{\partial t}{\partial q}$. 
Let
$$
T_n=[t_{0,n}, t_{1,n},t_{2,n},t_{3,n}, t_{4,n}, t_{5,n},t_{6,n}].
$$  
Comparing the coefficients of $q^0$ and $q^1$ in both sides of $\Ra$ we get:
$$
T_0=[\frac{1}{5},-25,-35,-6,0,-1, -15],\ 
$$
$$T_1=[24, -2250,-5350,-355, 1,1875,4675].
$$
Comparing the coefficients of $q^n, n\geq 2$ we get a recursion of the following type:
$$
(A_0+5nI_{7\times 7})T_n^\tr=\hbox {A function of the entries of } T_0,\ T_1, \ldots, T_{n-1},
$$
where
$$
A_0=[\frac{\partial (t_5\Ra_{i})}{\partial t_j}]_{i,j=0,1,\ldots,6} \hbox{  evaluated at } t=T_0,\ \ 
\Ra=\sum_{i=0}^6\Ra_i\frac{\partial }{\partial t_i}.
$$ 
The matrix $A_0+5nI_{7\times 7},\ n\geq 2$ is invertible and so we get a recursion in $T_n$.

\def\cprime{$'$} \def\cprime{$'$} \def\cprime{$'$}



\end{document}